\documentclass[11pt]{article}
\usepackage{amsmath,amsthm,amssymb}
\usepackage{enumerate}
\usepackage{graphicx}
\usepackage{titlesec}
\usepackage{empheq}
\usepackage{color}
\usepackage{xfrac,bigints}

\newtheorem{Th}{Theorem}[section]
\newtheorem{Lem}[Th]{Lemma} 
 
\newtheorem{Cor}[Th]{Corollary}

\newtheorem{Rem}[Th]{Remark}

\titleformat*{\section}{\normalsize\bfseries}
\titleformat*{\subsection}{\normalsize\bfseries}

\newcommand{\argmin}{\mathop{\rm arg~min}\limits}

\def\N{{\mathbb N}} 
\def\R{{\mathbb R}}

\def\setm{{\hspace{-0.3pt}\setminus \hspace{-0.3pt}}}
\def\sgn{\mbox{\rm sgn}}
\def\vep{\varepsilon}

\DeclareMathOperator*{\esssup}{ess\,sup}
\DeclareMathOperator*{\essinf}{ess\,inf}

\makeatletter 
\long\def\@makefntext#1{\parindent 1em\noindent 
\@hangfrom{\hbox to 1.8em{\hss$^{\@thefnmark}$}}#1}
\makeatother

\begin{document}

\begin{center}
  {\Large\bf
    Two-step minimization approach\\
    to an $L^\infty$-constrained variational problem\\
with a generalized potential
}
\end{center}
~\\
\begin{center}
  Vina Apriliani$^{\mbox{\scriptsize a,b}}$,
  Masato Kimura$^{\mbox{\scriptsize c}}$,
  Hiroshi Ohtsuka$^{\mbox{\scriptsize c}}$\\~\\
\begin{tabular}{l}
  {\scriptsize a)} Graduate School of Natural Science and Technology, Kanazawa University\\
 {\scriptsize b)} Universitas Islam Negeri Ar-Raniry\\
{\scriptsize c)} Faculty of Mathematics and Physics, Kanazawa University
\end{tabular}
\end{center}

\renewcommand{\thefootnote}{\fnsymbol{footnote}}
\footnote[0]{E-mail: vina.apriliani@ar-raniry.ac.id (VA), mkimura@se.kanazawa-u.ac.jp (MK), ohtsuka@se.kanazawa-u.ac.jp (HO)}
\renewcommand{\thefootnote}{\arabic{footnote}}

\vspace{0.5cm}
\begin{abstract}
We study a variational problem on $H^1(\R)$ under an $L^\infty$-constraint
related to Sobolev-type inequalities
for a class of generalized potentials, including $L^p$-potentials, non-positive potentials,
and signed Radon measures.
We establish various essential tools for this variational problem,
including the decomposition principle, the comparison principle,
and the perturbation theorem, which are the basis of
the two-step minimization method. As for their applications,
we present precise results for minimizers of minimization problems,
such as the
study of potentials of Dirac's delta measure type
and the analysis of trapped modes in potential wells.
\end{abstract}


\section{Introduction}\label{intro}
\setcounter{equation}{0}
In our previous paper \cite{AKOp1}, we considered the following variational problem
for a positive bounded potential $V\in L^\infty(\R)$ with $\essinf_{x\in\R} V(x) >0$:
\begin{align}\label{mmV}
m(V):=\min_{u\in H^1(\R),\|u\|_\infty = 1}
\int_{\R} (|u'(x)|^2+V(x)|u(x)|^2)\,dx,
\end{align}
where $\|u\|_\infty:=\|u\|_{L^\infty(\R)}$.
The above variational problem \eqref{mmV} is closely related
to the problem of finding the best constant for the following Sobolev-type inequality
for the $L^\infty$-norm:
\begin{equation}\label{So-ineq-inhom}
\|u\|_\infty\leq C\|u\|_V,
\end{equation}
where
\[
\|u\|_V^2
:=\int_{\R} (|u'(x)|^2+V(x)|u(x)|^2)\,dx.
\]
The best constant of the Sobolev-type inequality \eqref{So-ineq-inhom}
is given by $m(V)^{-1/2}$.
The Sobolev-type inequality has been studied in various settings \cite{CTZ19,CM16,Kam08,WKNTY08}.
We proposed a two-step minimization approach to the variational problem \eqref{mmV} in \cite{AKOp1}.
This approach enabled us
to evaluate the best constant for inhomogeneous positive bounded potentials precisely.

The aim of this paper is to consider the variational problem \eqref{mmV}
with a more general class of potentials $V$,
including unbounded potentials, non-positive potentials,
and the Dirac delta measures,
and to extend the two-step minimization approach to them.
We provide various essential tools for this variational problem,
such as the decomposition principle that provides the basis
for the two-step minimization method,
the comparison principle,
the perturbation theorem, 
the existence and some properties of the minimizer, the continuity of the minimum value
in the first minimization step, etc.

As applications of those tools, we consider two specific potential cases.
For a potential that contains the Dirac delta measure,
applying our  comparison principle, we give an alternative proof
for the best constant of the Sobolev-type inequality
\[
\|u\|_\infty \le C
\left( \int_{\R} \left( |u'(x)|^2
+\alpha |u(x)|^2\right) dx+\beta |u(0)|^2\right)^\frac{1}{2},
\]
by \cite{Kam08} in Theorem~\ref{kam2}.
The other application is devoted to the case for a potential well.
By using the precise property of the minimizer, we will provide
sufficient conditions for the trapped mode in terms of the
depth and width of the potential well.

The structure of this paper is as follows.
Section~\ref{sec:3.1} introduces a class of generalized potentials and
demonstrates the decomposition principle underlying
the two-step minimization approach.
In Section~\ref{sec:2.1}, we survey the results of \cite{AKOp1} for positive bounded potentials,
which are necessary for the discussion of this paper.
Section~\ref{sec:cp-pt} describes a comparison principle and a perturbation theorem of $m(V)$
for generalized potentials.
As their application, we give an alternative proof 
for the potential given by a constant plus Dirac delta measure.
In Section~\ref{sec:pm},
we study 
the existence of the minimizer in the first minimization step
and also the continuity of the minimum value.
We finally provide a sufficient condition for the trapped mode
in the potential well.

\section{Generalized potential and decomposition principle}\label{sec:3.1}
In this paper, we deal with Sobolev-type inequalities
related to the Schr\"{o}dinger-type operator $-\frac{d^2}{dx^2}+V$
with a potential term $V$.
We begin this section by introducing a general class of the potential $V$.
We define
\[
X:=\{ V: H^1(\R)\times H^1(\R)\to \R;~
\mbox{$V$: bounded symmetric bilinear map}\}.
\]
Then, it is known that $X$ is a Banach space over $\R$ with the norm:
\[
\|V\|_X:=\sup_{u,v\in H^1(\R),u,v\not{\equiv}0}
\frac{|V(u,v)|}{\|u\|_{H^1(\R)}\|v\|_{H^1(\R)}}.
\]

Without loss of generality, 
we suppose that $u\in H^1(\R)$ (or more generally, $u\in W^{1,p}_{loc}(\R)$
for $p\in [1,\infty]$) always satisfies $u\in C^0(\R)$, since
an element of the function space $W^{1,p}_{loc}(\R)$ has a continuous representation
(Theorem~8.8 of \cite{Bre11}). We also remark that
$u\in H^1(\R)$ satisfies $u\in L^\infty(\R)$ and
$\lim_{|x|\to\infty}u(x)=0$ (Theorem~8.8 and Corollary~8.9 of \cite{Bre11}).

We remark that $L^p(\R)$ ($1\le p\le \infty$) is continuously embedded
in $X$ by identifying $V\in L^p(\R)$ with the following $\tilde{V}\in X$:
\[
\tilde{V}(u,v):=\int_\R V(x)u(x)v(x)\,dx
\quad
(u,v\in H^1(\R)),
\]
where $Vuv\in L^1(\R)$ is clear from $uv\in L^1(\R)\cap L^\infty(\R)\subset L^q(\R)$
for any $q\in [1,\infty]$.

We denote the set of signed Radon measure $V$ on $\R$ with finite total variation $|V|_{TV}:=|V|(\R)<\infty$
by ${\cal M}_1(\R)$. For example, the Dirac measure $\delta_a$ belongs to ${\cal M}_1(\R)$.
It is defined by $\delta_a(A)=1$ if $a\in A$,
and $\delta_a(A)=0$ if $a\not\in A$ for $a\in \R$ and $A\in {\cal B}(\R)$,
where ${\cal B}(\R)$ is the set of Borel sets on $\R$.
We also remark that $L^1(\R)\subset {\cal M}_1(\R)$ and
$|V|_{TV}=\| V\|_{L^1(\R)}$ holds for $V\in L^1(\R)$.

Then, $V\in {\cal M}_1(\R)$ 
is also considered as $V\in X$
by identifying $V$ with the following $\tilde{V}\in X$:
\[
\tilde{V}(u,v):=\int_\R uv\,dV
\quad
(u,v\in H^1(\R)).
\]
It is well-defined since $uv\in C^0(\R)\cap L^\infty(\R)$
and $\|\tilde{V}\|_X\le 2 |V|_{TV}$ holds since
\begin{align*}
  \left| \tilde{V}(u,v) \right|
  \le \| u\|_\infty \| v\|_\infty |V|_{TV}\le 2\| u\|_{H^1(\R)} \| v\|_{H^1(\R)} |V|_{TV},
\end{align*}
where the last inequality follows from $\|u\|_\infty \le \sqrt{2}\| u\|_{H^1(\R)}$
(see p.213 of \cite{Bre11}). We note that the best constant of this inequality
is $\|u\|_\infty \le \frac{1}{\sqrt{2}}\| u\|_{H^1(\R)}$ (see \cite{AKOp1}).

For $V\in X$, we define 
\begin{align*}
  I(u; V):=\| u'\|_{L^2(\R)}^2+V(u,u)\quad (u\in H^1(\R)),
\end{align*}
and define a Rayleigh-type quotient:
\begin{align*}
  R(u; V):=\frac{I(u;V)}{\|u\|_{\infty}^2}
\quad (u\in H^1(\R)\setminus \{0\}). 
\end{align*}
We also define $m(V)\in [-\infty,\infty)$ and $M(V)$ by
\begin{align}
&m(V):=\inf_{u\in H^1(\R),u\not{\equiv}0} R(u; V),\label{mV}\\
&M(V):=\left\{u\in H^1(\R)\setminus \{ 0\};~m(V)=
R(u; V)\notag
\right\}.
\end{align}
Then, if and only if $m(V)>0$,
the following Sobolev-type inequality holds:
\begin{align}\label{sti}
  ~^\exists C>0~~\text{s.t}~~\| u\|_\infty \le C \,I(u;V)^{\frac{1}{2}}\quad
  (^\forall u\in H^1(\R)).
\end{align}
In this case,  $C=m(V)^{-1/2}$
gives the best constant of Sobolev-type inequality \eqref{sti}.

In this paper, we often consider the class of the generalized potentials
$L^\infty(\R)+{\cal M}_1(\R)\subset X$, where
\[
L^\infty(\R)+{\cal M}_1(\R)=\{ V=V_0+V_1\in X;~V_0\in L^\infty(\R),~V_1\in {\cal M}_1(\R)\}.
\]
We remark that $L^p(\R)\subset L^\infty(\R)+{\cal M}_1(\R)$ holds
for any $p\in [1,\infty ]$.
Indeed, it is trivial if $p=\infty$, and if $p\in [1,\infty)$,
let $V\in L^p(\R)$ and set
$A:=\{x\in\R;~|V(x)|>1\}$, $V_0(x):=(1-\chi_A(x))V(x)$, and $V_1(x):=\chi_A(x)V(x)$,
where $\chi_A$ is the indicator function of $A$.
Then, $\|V_0\|_\infty\le 1$ and
\[
\|V_1\|_{L^1(\R)}=\int_A |V(x)|\,dx\le \int_A |V(x)|^p\,dx \le \|V\|_{L^p(\R)}^p<\infty.
\]
Therefore, we obtain
\[
V=V_0+V_1\in L^\infty(\R)+L^1(\R)\subset L^\infty(\R)+{\cal M}_1(\R),
\]
since $V_0\in L^\infty(\R)$ and $V_1\in L^1(\R)\subset {\cal M}_1(\R)$.

In \cite{AKOp1},
the authors proved the following decomposition principle of the minimization problem
\eqref{mV} for general non-constant bounded positive potentials $V\in L^\infty(\R)$,
and established the two-step minimization approach to study
the precise properties of the minimizer for the Sobolev-type inequality.
We aim to extend the two-step minimization approach
to the case of the generalized potential $V\in X$ in this paper.

For $a\in \R$, we set
\[
K_a:=\{ u\in H^1(\R);~u(a)=\|u\|_\infty =1\},
\]
and define
\[
F(a;V):=\inf_{u\in K_a} I(u;V).
\]
We remark that $K_a$ is a closed convex set in $H^1(\R)$, which implies
that $K_a$ is weakly closed in $H^1(\R)$.
\begin{Th}[decomposition principle]\label{decomposition}
Let $V\in X$ and set $m(V)$ as \eqref{mV}. Then, we have
\begin{align}\label{infF}
m(V)
=\inf_{a\in\R}F(a;V).
\end{align}
\end{Th}
\proof
We first remark that $F(a;V)$ and $m(V)$
can have their values in $[-\infty,\infty)$.
For $a\in\R$, there exists $\{ u_{a,n}\}_{n=1}^\infty \subset K_a$
such that $\lim_{n\to\infty}I(u_{a,n};V)=F(a;V)$. Since
\[
m(V)\le R(u_{a,n};V)=I(u_{a,n};V),
\]
it follows that $m(V)\le F(a;V)$.
Let us define $\tilde{m}(V):=\inf_{a\in\R}F(a;V)$.
Then, taking the infimum concerning $a$ in $m(V)\le F(a;V)$, we obtain $m(V)\le \tilde{m}(V)$.

Let $\{u_n\}_{n=1}^\infty$ be a minimizing sequence attaining the infimum of \eqref{mV}.
Choosing $a_n\in\R$ as $|u_n(a_n)|=\|u_n\|_\infty >0$,
we define $v_n:=u_n(a_n)^{-1}u_n\in K_{a_n}$.
Then we have
\[
m(V)=\lim_{n\to\infty}R(u_n;V)=\lim_{n\to\infty}I(v_n;V),
\]
and $\tilde{m}(V)\le m(V)$ follows from $\tilde{m}(V)\le F(a_n;V)\le I(v_n;V)$
as $n\to\infty$. Hence, we obtain $\tilde{m}(V)=m(V)$.
\qed

\section{Bounded positive potentials}\label{sec:2.1}
\setcounter{equation}{0}
We briefly summarize the results obtained in
\cite{AKOp1} for the case of bounded positive potentials.
In this section, we suppose 
\begin{align}\label{bddV}
V\in L^\infty(\R),\quad
0<v_0:=\essinf_{x\in\R}V(x),\quad v_1:=\esssup_{x\in\R}V(x).
\end{align}
Then, for $u,~v\in H^1(\R)$, we define
\begin{align*}
  &(u,v)_V:=\int_\R \left(u'(x)v'(x)+V(x)u(x)v(x)\right)\,dx,\quad \| u\|_V:=(u,u)_V^{\frac{1}{2}}.
\end{align*}
We remark that $(u,v)_V$ defines an inner product on $H^1(\R)$.
The corresponding norm $\|u\|_V$ is equivalent to the norm
of $H^1(\R)$ and it satisfies $I(u;V)=\|u\|_V^2$.

We consider the first minimization step:
\begin{align}\label{Fa}
F(a;V)= \inf_{u\in K_a}\| u\|_V^2.
\end{align}

\begin{Th}[\cite{AKOp1}]\label{th1}
We suppose the condition \eqref{bddV}
and fix $a\in\R$. There exists a unique minimizer $u_a\in K_a$ to \eqref{Fa}, that is,
\begin{align}
  &u_a=\argmin\limits_{u\in K_a}\|u\|_V^2,\label{uaua}\\
  &F(a;V)= \min_{u\in K_a}\| u\|_V^2=\| u_a\|_V^2\notag ,
\end{align}
and it satisfies the following properties:
\begin{align}
&u_a\in W^{2,\infty}(\R\setm\{a\})\quad \mbox{and}\quad
u_a''(x)=V(x)u_a(x)~(\mbox{a.e.}~ x\in\R\setm \{a\}).\label{p1}\\
&  
e^{-\sqrt{v_1}|x-a|}\le u_a(x) \le e^{-\sqrt{v_0}|x-a|}\quad(x\in\R),\label{p2}\\
&
\frac{v_0}{\sqrt{v_1}}e^{-\sqrt{v_1}|x-a|}\le \sgn (a-x)  u_a'(x)
\le \frac{v_1}{\sqrt{v_0}}e^{-\sqrt{v_0}|x-a|}\quad(x\in\R\setm \{a\}).\label{p3}
\end{align}
\end{Th}

\begin{Th}[\cite{AKOp1}]\label{nondecreasingV}
  We assume \eqref{bddV} and suppose that $V$ is a non-decreasing function.
Then, it holds that
$m(V)=\lim_{a\to -\infty}F(a;V) =2\sqrt{v_0}$.
Furthermore, if $v_0<v_1$, then $F$ is a strictly increasing function
and $M(V)=\emptyset$. If\, $V$ is constant,
then it holds that
$m(V)=2\sqrt{V}$ and 
$M(V)=\{
cu_a;~c\in\R\setm\{0\},~a\in \R
\}$,
where $u_a(x)=e^{-\sqrt{V}|x-a|}$.
\end{Th}

\section{Comparison principle and perturbation theorem}\label{sec:cp-pt}
From this section onwards, 
we consider the generalized potentials.
We consider the following comparison principle of $m(V)$.
\begin{Th}[comparison principle of $m(V)$]\label{cpmV}
  We suppose that $V_1,V_2\in X$ and that $m(V_1)\ne -\infty$ or $m(V_2)\ne -\infty$.
  Then we have
\begin{align}\label{cpsbc}
  \inf \limits_{u\in H^1(\R),u\ne 0} \frac{(V_1-V_2)(u,u)}{\|u\|_\infty^2}
  \le m(V_1)-m(V_2)\le
  \sup \limits_{u\in H^1(\R),u\ne 0} \frac{(V_1-V_2)(u,u)}{\|u\|_\infty^2}.
\end{align}
\end{Th}
\proof
For $V_1$, let $\{u_n\}_{n\in\N}\subset H^1(\R)$ be a minimizing sequence
to the infimum $m(V_1)$ and suppose $\|u_n\|_\infty =1$.
Then, it satisfies
$\lim_{n\to \infty}I(u_n;V_1)=m(V_1)$ and 
\begin{align*}
  I(u_n;V_1)-m(V_2)\ge I(u_n;V_1)-I(u_n;V_2)
  =(V_1-V_2)(u_n,u_n)
  \ge \inf \limits_{u\in H^1(\R),u\ne 0} \frac{(V_1-V_2)(u,u)}{\|u\|_\infty^2}.
\end{align*}
Taking the limit as $n\to\infty$, we obtain the first inequality of \eqref{cpsbc}.
By exchanging $V_1$ and $V_2$, we derive the second inequality as
\begin{align*}
m(V_1)-m(V_2)
\le -\left(\inf \limits_{u\in H^1(\R),u\ne 0} \frac{(V_2-V_1)(u,u)}{\|u\|_\infty^2}\right)
=
\sup \limits_{u\in H^1(\R),u\ne 0} \frac{(V_1-V_2)(u,u)}{\|u\|_\infty^2}.
\end{align*}
\qed

\begin{Cor}\label{newcor}
  Under the condition of Theorem~\ref{cpmV}, if $(V_1-V_2)(u,u)\ge 0$
  for $u\in H^1(\R)$, then $m(V_2)\le m(V_1)$ holds.
\end{Cor}

Also, from Theorem~\ref{cpmV}, we immediately have the following theorem.
\begin{Th}\label{cor4.2}
 We suppose that $V\in X$ and $m(V)\ne -\infty$.
 If $\mu \in {\cal M}_1(\R )$, then
\begin{align}\label{mu+-}
-|\mu_- |_{TV}\le m(V+\mu)-m(V) \le |\mu_+ |_{TV},
\end{align}
where $\mu_+$ and $\mu_-$ are the positve and negative parts of the Radon measure $\mu$.
In particular, we have
\[
|m(V+\mu)-m(V)| \le \max \big(|\mu_- |_{TV},|\mu_+ |_{TV}\big)\le |\mu |_{TV}.
\]
\end{Th}
\proof
We apply Theorem~\ref{cpmV} with $V_1=V+\mu$ and $V_2=V$.
Then, we have
\begin{align*}
  \inf \limits_{u\in H^1(\R),u\ne 0} \frac{\mu (u,u)}{\|u\|_\infty^2}
  \le m(V_1)-m(V_2)\le
  \sup \limits_{u\in H^1(\R),u\ne 0} \frac{\mu (u,u)}{\|u\|_\infty^2}.
\end{align*}
Since $\mu (u,u)=\mu_+(u,u)-\mu_-(u,u)$, 
paying attention to the following inequalities
\[
-|\mu_-|_{TV}\|u\|_\infty^2\le -\mu_- (u,u)
\le \mu(u,u)\le \mu_+(u,u) \le |\mu_+|_{TV}\|u\|_\infty^2,
\]
we obtain \eqref{mu+-}. The last inequality also follows from
$|\mu |_{TV}=|\mu_+|_{TV}+|\mu_-|_{TV}$.
\qed

\begin{Th}\label{limv}
  Let $V\in L^\infty(\R)$ with $\essinf_{x\in\R}V(x)>0$.
  We suppose
\begin{align}\label{limFa}
  m(V)=\lim_{a\to\infty}F(a;V),\quad \mbox{or}\quad m(V)=\lim_{a\to -\infty}F(a;V).
\end{align}
Let $p\in [1,\infty)$.
If $\mu\in L^p(\R)+{\cal M}_1(\R)\subset X$ is nonnegative,
i.e., $\mu(u,u)\ge 0$ for $u\in H^1(\R)$,
  then $m(V+\mu)=m(V)$ and $M(V+\mu)\subset M(V)$ hold.
\end{Th}
\proof
Choosing $V_1=V+\mu$ and $V_2=V$ in Theorem~\ref{cpmV}, we have
\begin{align}\label{q1}
m(V+\mu)-m(V)\ge \inf \limits_{u\in H^1(\R),u\ne 0} \frac{\mu (u,u)}{\|u\|_\infty^2}\ge 0.
\end{align}

We define $u_a(x)$ as in Theorem~\ref{th1}. Then, we have
\[
m(V+\mu)\le I(u_a;V+\mu)=I(u_a;V)+\mu(u_a,u_a)=F(a;V)+\mu(u_a,u_a).
\]
From the assumption \eqref{limFa}, taking the limit $a\to \infty$ or $a\to -\infty$,
we obtain
\begin{align}\label{q2}
m(V+\mu)\le m(V)+\lim_{a\to \pm\infty}\mu(u_a,u_a)=m(V),
\end{align}
where the last equality holds as follows.

Let $\mu =\mu_0+\mu_1$ with $\mu_0\in L^p(\R)$ and $\mu_1\in {\cal M}_1(\R)$.
Since $L^1(\R)\subset {\cal M}_1(\R)$, we assume $p\in (1,\infty)$ without loss of generality
and define $q\in (1,\infty)$ as $p^{-1}+q^{-1}=1$.
From the estimate \eqref{p2}, we have $|u_a(x)|\le e^{-\sqrt{v_0}|x-a|}$.
For $R>0$, we define $I_R:=[-R,R]$ and $J_R:=\R\setminus I_R$,
and suppose $a\in J_R$.
Then,
$|u_a(x)|\le e^{-\sqrt{v_0}(|a|-R)}$ holds for $x\in I_R$.
Hence, we have
\begin{align*}
  \mu_0(u_a,u_a)&=\int_{I_R}\mu_0(x)|u_a(x)|^2\,dx
  + \int_{J_R}\mu_0(x)|u_a(x)|^2\,dx\\
  &\le e^{-2\sqrt{v_0}(|a|-R)}\| \mu_0\|_{L^1(I_R)}
  + \|\mu_0\|_{L^p(J_R)}\|u_a^2\|_{L^q(\R)}.
\end{align*}
Noting that 
\begin{align*}
  \|u_a^2\|_{L^q(\R)}\le \left(\int_{\R}e^{-2q\sqrt{v_0}|x-a|}\,dx\right)^{\frac{1}{q}}
  =\left(q\sqrt{v_0}\right)^{-\frac{1}{q}},
\end{align*}
for an arbitrary $\vep >0$, there exists $R>0$ such that
$\|\mu_0\|_{L^p(J_R)}\|u_a^2\|_{L^q(\R)}\le \vep$. Then, there exists $\tilde{R}>R$ such that
$e^{-2\sqrt{v_0}(|a|-R)}\| \mu_0\|_{L^1(I_R)}\le\vep$ holds for $|a|>\tilde{R}$.
It implies that $\lim_{|a|\to\infty}\mu_0(u_a,u_a)=0$.

Similarly, we have
\begin{align*}
  \mu_1(u_a,u_a)&=\int_{I_R}|u_a|^2\,d\mu_1
  + \int_{J_R}|u_a|^2\,d\mu_1\le e^{-2\sqrt{v_0}(|a|-R)}|\mu_1|(I_R)
+|\mu_1|(J_R).
\end{align*}
For an arbitrary $\vep >0$, there exists $R>0$ such that
$|\mu_1|(J_R)\le \vep$. Then, there exists $\tilde{R}>R$ such that
$e^{-2\sqrt{v_0}(|a|-R)}|\mu_1|(I_R)\le\vep$ holds for $|a|>\tilde{R}$.
It implies that $\lim_{|a|\to\infty}\mu_1(u_a,u_a)=0$.

Hence, we conclude $m(V+\mu)=m(V)$ from \eqref{q1} and \eqref{q2}.
Moreover, for $u\in M(V+\mu)$, since we have
\[
I(u;V)\le I(u;V)+\mu(u,u)=I(u;V+\mu)=m(V+\mu)=m(V)\le I(u;V),
\]
$I(u;V)=m(V)$ follows and it implies the inclusion $M(V+\mu)\subset M(V)$.
\qed

Using Theorem~\ref{limv},
we can give an alternative proof for the following result.
We define $\delta_0\in {\cal M}_1(\R)\subset X$ by $\delta_0(u,v):=u(0)v(0)$.
\begin{Th}[Kametaka et al. \cite{Kam08}]\label{kam2}
Let $\alpha  > 0$ and $\beta \in\R$.
Then $m(\alpha +\beta \delta_0)=2\sqrt{\alpha}-\beta_-$ holds,
where $\beta_-=\max (-\beta,0)$.
In particular, if $2\sqrt{\alpha}+\beta >0$, then the Sobolev-type inequality
\[
\|u\|_\infty \le C
\left( \int_{\R} \left( |u'(x)|^2
+\alpha |u(x)|^2\right) dx+\beta |u(0)|^2\right)^\frac{1}{2}
\]
holds and its best constant is given by $C=(2\sqrt{\alpha}-\beta_-)^{-\frac{1}{2}}$.
\end{Th}
\proof
For the case of $\beta \ge 0$, $m(\alpha+\beta\delta_0)=m(\alpha)=2\sqrt{\alpha}
=2\sqrt{\alpha}-\beta_-$
holds from Theorems~\ref{nondecreasingV} and \ref{limv}, since $V=\alpha$ satisfies the condition \eqref{limFa}.

If $\beta<0$, from Theorem~\ref{cpmV} with $V_1=\alpha+\beta\delta_0$ and $V_2=\alpha$,
we have
\begin{align}\label{e1}
m(\alpha+\beta\delta_0)-m(\alpha)\ge 
\inf \limits_{u\in H^1(\R),u\ne 0} \frac{\beta |u(0)|^2}{\|u\|_\infty^2}
=\beta \left(\sup\limits_{u\in H^1(\R),u\ne 0} \frac{|u(0)|^2}{\|u\|_\infty^2}\right)
=\beta .
\end{align}
On the other hand, setting $u_0(x)=e^{-\sqrt{\alpha}|x|}$, we also have
\begin{align}\label{e2}
  m(\alpha+\beta\delta_0)\le I(u_0;\alpha+\beta\delta_0)
  =I(u_0;\alpha)+\beta |u_0(0)|^2=m(\alpha)+\beta .
\end{align}
Hence, from \eqref{e1} and \eqref{e2}, we obtain
$m(\alpha+\beta\delta_0)=m(\alpha)+\beta=2\sqrt{\alpha}-\beta_-$.
\qed

\section{Properties of the minimizers}\label{sec:pm}
\setcounter{equation}{0}
In the following discussion, we will often
make the following assumption on the generalized potential $V\in X$:
\begin{align}\label{V01}
  V=V_0+V_1,~~V_0\in L^\infty(\R),~~\essinf\, V_0 >0,~~V_1\in {\cal M}_1(\R).
\end{align}
We note that the next identity holds under the assumption \eqref{V01},
\begin{align}\label{IV01}
I(u;V)=I(u;V_0)+V_1(u,u)\quad (u\in H^1(\R)).
\end{align}

\begin{Lem}\label{newlem}
We suppose that $V\in X$ satisfies \eqref{V01}.
\begin{enumerate}
\item
  There exist $C_1$, $C_2>0$ such that
  the following inequality holds for $u\in H^1(\R)$:
\begin{align}\label{L1}
  \| u\|_{H^1(\R)}^2\le C_1I(u;V)+C_2\| u\|_\infty^2.
\end{align}
\item
  $I(\cdot ;V)$ is weakly lower semi-continuous in $H^1(\R)$.
\end{enumerate}
\end{Lem}
\proof
We set $\alpha :=\essinf\, V_0 >0$.
For $u\in H^1(\R)$, using \eqref{IV01}, we have
\begin{align*}
\min (1,\alpha)\| u\|_{H^1(\R)}^2&\le \int_\R \left(|u'(x)|^2+\alpha|u(x)|^2\right)\,dx\notag\\
&\le I(u;V_0)\notag\\
  &= I(u;V)-V_1(u,u)\notag\\
  &\le I(u;V)+|V_1|_{TV}\|u\|_\infty^2.
\end{align*}
Hence, we have the first assertion \eqref{L1} by setting
\begin{align*}
  C_1:=\frac{1}{\min (1,\alpha)},\quad
  C_2:=\frac{|V_1|_{TV}}{\min (1,\alpha)}.
\end{align*}

Let us suppose that $\{u_n\}_{n\in\N}\subset H^1(\R)$ and $u\in H^1(\R)$
satisfy
$u_n\to u$ weakly in $H^1(\R)$ as $n\to\infty$.
Since $\essinf\, V_0 >0$, the inner product
$(\cdot,\cdot)_{V_0}$ gives
an equivalent topology on $H^1(\R)$ and $I(u;V_0)=(u,u)_{V_0}$ holds.
Therefore, it follows that
$I(\cdot ;V_0)$ is weakly lower semi-continuous in $H^1(\R)$, i.e.,
it holds that
\begin{align}\label{limV0}
I(u;V_0)\le \liminf_{n\to\infty}I(u_n;V_0).
\end{align}

Paying attention to the identity \eqref{IV01},
it is sufficient to show that
\begin{align}\label{V1eq}
\lim_{n\to\infty} V_1(u_n,u_n)=V_1(u,u),
\end{align}
to prove the second assertion of the lemma.

For any $\vep >0$,
from $|V_1|_{TV}<\infty$,
there exists $R>0$ such that
\begin{align}\label{vep1}
\int_{\{|x|\ge R\}} \,d|V_1|\le \vep.
\end{align}
Since $\{u_n\}_{n\in\N}$ is boudned in $H^1(\R)$,
from the Rellich-Kondrachov theorem for the compact embedding
$H^1(-R,R)\subset C^0([-R,R])$,
there exists a subsequence which uniformly convergent on $[-R,R]$.
However, the limit function of the uniform convergence 
coincides with $u$. As a result, the whole sequence $\{u_n\}_{n\in\N}$
converges uniformly to $u$ on $[-R,R]$.
Hence,
there exists $N\in\N$ such that
\begin{align}\label{vep2}
  \|u_n-u\|_{L^\infty(-R,R)}<\vep \quad (n\ge N).
\end{align}

We set
$B:=\sup_{n\in\N}\|u_n\|_\infty <\infty$.
It implies $\|u\|_\infty\le B$.
Thus, from \eqref{vep1} and \eqref{vep2}, we obtain
\begin{align*}
  \big| V_1(u_n,u_n)-V_1(u,u)\big|
  &= \left| \int_\R \left( |u_n|^2-|u|^2\right)\,dV_1\right|\\
&  \le \int_{\{|x|<R\}} \left|u_n+u\right|\,\left| u_n-u\right|\,d|V_1|
    +\int_{\{|x|\ge R\}} \left( |u_n|^2+|u|^2\right)\,d|V_1|\\
&    \le  2B| V_1|_{TV}\|u_n-u\|_{L^\infty(-R,R)}+2B^2\int_{\{|x|\ge R\}} \,d|V_1|\\
&    \le 2B\left(| V_1|_{TV}+B\right)\vep.
\end{align*}
From this estimate, we obtain \eqref{V1eq} and
\begin{align*}
  I(u;V)=I(u;V_0)+V_1(u,u)
  \le \liminf_{n\to\infty}I(u_n;V_0)+\lim_{n\to\infty} V_1(u_n,u_n)
  =\liminf_{n\to\infty}I(u_n;V).
\end{align*}
\qed

\begin{Th}\label{fms1}
We suppose \eqref{V01}.
Then, for each $a\in\R$, there exists $u_a\in K_a$ such that
\begin{align}\label{Iua}
I(u_a;V)=\min_{u\in K_a} I(u;V),
\end{align}
i.e., $F(a;V)$ is attained as $F(a;V)=I(u_a;V)$.
\end{Th}
\proof
For any fixed $a\in\R$,
let $\{ u_{a,n}\}_{n\in\N} \subset K_a$
be a minimizing sequence for $I(\cdot;V)$ in $K_a$, i.e.,
$\lim_{n\to\infty}I(u_{a,n};V)=F(a;V)$.
From \eqref{L1}, $\{ u_{a,n}\}_{n\in\N}$ is bounded in $H^1(\R)$.
Replacing $\{ u_{a,n}\}_{n\in\N}$ by a subsequence if necessary,
there exists $u_a\in H^1(\R)$ such that $u_{a,n}$ weakly converges
to $u_a$ in $H^1(\R)$ as $n\to \infty$.
Since $u_{a,n}\in K_a$ and $K_a$ is weakly closed, $u_a\in K_a$ holds.
Hence, we obtain \eqref{Iua} as 
$I(u_a;V)\le \lim_{n\to\infty} I(u_{a,n};V)=F(a;V)$
from
the second claim of Lemma~\ref{newlem}.
\qed

\begin{Th}\label{fms2}
We suppose that \eqref{V01} holds with $V_1\in L^1(\R)$,
and that $u_a\in K_a$ satisfies \eqref{Iua}.
Then, $u_a(x)>0$ holds for $x\in\R$, and, setting $J:=\{x\in\R;~u_a(x)<1\}$, 
$u_a\in W^{2,1}(J)$ and $u_a''(x)=V(x)u_a(x)$ hold for a.e. $x\in J$.  
\end{Th}
\proof
For $u_a\in K_a$ which satisfies \eqref{Iua}, we define $J:=\{x\in\R;~|u_a(x)|<1\}$,
which is an open set in $\R$ since $u_a$ is continuous.
Then, for any $\varphi\in C^\infty_0(J\,)$, there exists $\tau >0$ such that
$u_a+t\varphi\in K_a$ for $t\in (-\tau,\tau)$.
Since $I(u_a+t\varphi;V)$ has a local minimum at $t=0$, we obtain
\begin{align*}
  0=\frac{d}{dt}I(u_a+t\varphi;V)\big|_{t=0}
  =2\int_{J}\left(u_a'(x)\varphi'(x)+V(x)u_a(x)\varphi(x)\right)\,dx,
\end{align*}
which implies that $u_a\in W^{2,1}(J\,)$
and $u_a''(x)=V(x)u_a(x)$ holds for a.e. $x\in J$.

We set $\bar{u}_a(x):=|u_a(x)|$. Then, since $\bar{u}_a\in K_a$ and
$I(\bar{u}_a;V)=I(u_a;V)=F(a;V)$, we have that
$\bar{u}_a\in W^{2,1}(J\,)$ and
\begin{align}\label{baru}
  \bar{u}_a''(x)=V(x)\bar{u}_a(x)\quad
  (\mbox{a.e. } x\in J\,).
\end{align}

Let $J_0$ be an open component of $J$.
Then,
\begin{align}\label{JJ}
  \bar{u}_a(x)=1\quad (x\in \overline{J_0}\setminus J_0\neq \emptyset),
\end{align}
from the definition of $J$.

We assume that $u_a(x_0)=0$ holds at some $x_0\in J_0$.
Then $\bar{u}_a(x_0)=\bar{u}_a'(x_0)=0$ holds,
since $\bar{u}_a\in C^1(J)$ and $\bar{u}_a(x)\ge 0$.
From \eqref{baru}, for $x\in J_0\cap [x_0,\infty)$, we have
\begin{align*}
  \bar{u}_a(x)&=\int_{x_0}^x\bar{u}_a'(s)\,ds
  \le  \int_{x_0}^x|\bar{u}_a'(s)|\,ds,\\
  |\bar{u}_a'(x)|&=\left|\int_{x_0}^x\bar{u}_a''(s)\,ds\right|
  \le  \int_{x_0}^x|V(s)|\bar{u}_a(s)\,ds,
\end{align*}
Setting $v(x):=\bar{u}_a(x)+|\bar{u}_a'(x)|$,
we have
\begin{align}\label{vineq}
  v(x)\le \int_{x_0}^x(1+|V(s)|)v(s)\,ds
  \quad (x\in J_0\cap [x_0,\infty)).
\end{align}
Applying the Gronwall inequality to \eqref{vineq},
we obtain that $v(x)\le 0$ for $x\in J_0\cap [x_0,\infty)$.
Since we can similarly obtain $v(x)\le 0$ for $x\in J_0\cap (-\infty, x_0]$ too,
$v(x)=0$ holds for $x\in J_0$.
However, this contradicts \eqref{JJ}.
Hence, we conclude that $u_a(x)\neq 0$ for $x\in\R$.
It implies that $u_a(x)=\bar{u}_a(x)>0$ for $x\in\R$ and
$J=\{x\in\R;~u_a(x)<1\}$.
\qed

\begin{Lem}\label{Lebsgue}
  We suppose $V\in L^\infty(\R)+{\cal M}_1(\R)$ and $u\in H^1(\R)$.
  Then, it holds that
  \begin{align}\label{lemlem}
  \lim_{h\to 0}I(u(\cdot -h);V)=I(u;V).
  \end{align}
\end{Lem}
\proof
First, for $x,~h\in\R$ and $u\in H^1(\R)$, we remark that
\[
|u(x)-u(x-h)|=\left| \int_{x-h}^x u'(y)\,dy\right|
\le \left( \left| \int_{x-h}^x|u'(y)|^2\,dy\right|\right)^{\frac{1}{2}}
|h|^{\frac{1}{2}}
  \le
  \|u\|_{H^1(\R)}|h|^{\frac{1}{2}}.
    \]
We use the idea in the proof of
Proposition~4.2.6 of \cite{Will07}.
For $h\in \R$, we set
\[
f_h(x):=|u(x)|+|u(x-h)|-|u(x)-u(x-h)|\ge 0\quad (x\in\R).
\]
Since $u$ is continuous, $\lim_{h\to 0}f_h(x)=2|u(x)|$ holds for $x\in \R$.
Applying Fatou's lemma, we obtain
\begin{align*}
2\int_\R |u(x)|\,dx&\le \liminf_{h\to 0}\int_\R f_h(x)\,dx\\
&=\liminf_{h\to 0}\int_\R \big( |u(x)|+|u(x-h)|-|u(x)-u(x-h)|\big) \,dx\\
&=2\int_\R |u(x)|\,dx+\liminf_{h\to 0}\int_\R \big( -|u(x)-u(x-h)|\big) \,dx\\
&=2\int_\R |u(x)|\,dx-\limsup_{h\to 0}\int_\R |u(x)-u(x-h)| \,dx.
\end{align*}
This implies $\limsup_{h\to 0}\int_\R |u(x)-u(x-h)| \,dx \le 0$
and also
\[
\lim_{h\to 0}\int_\R |u(x)-u(x-h)| \,dx =0.
\]

We write $V=V_0+V_1$, where $V_0\in L^\infty(\R)$ and $V_1\in {\cal M}_1(\R)$.
Then we have
\begin{align*}
  I(u;V)-I(u(\cdot -h);V)
  &= V(u,u)-V(u(\cdot -h),u(\cdot -h))\\
  &=\int_\R V_0(x)\left( |u(x)|^2-|u(x-h)|^2\right)\,dx
  + \int_\R \left( |u|^2-|u(\cdot -h)|^2\right)\,dV_1.
\end{align*}
Hence, we obtain \eqref{lemlem} from 
\begin{align*}
  \big| I(u;V)-I(u(\cdot -h);V)\big|
  &\le 2\| V_0\|_\infty \|u\|_\infty  \int_\R |u(x)-u(x-h)|\,dx
  + 2\|u\|_\infty \int_\R |u-u(\cdot -h)|\,d|V_1|\\
  &\le 2\| V_0\|_\infty \|u\|_\infty  \int_\R |u(x)-u(x-h)|\,dx
  + 2\|u\|_\infty |V_1|_{TV}\|u\|_{H^1(\R)}|h|^{\frac{1}{2}}.
\end{align*}
\qed
\begin{Th}\label{FC0}
If\, $V\in X$ satisfies \eqref{V01}, then $F(\cdot;V)\in C^0(\R)$ holds.
\end{Th}
\proof
For $a\in\R$ and any convergent sequence $a_n\to a$ as $n\to\infty$,
from Theorem~\ref{fms1},
there exists $u_a\in K_a$ and $u_{a_n}\in K_{a_n}$ such that
$F(a;V)=I(u_a;V)$ and $F(a_n;V)=I(u_{a_n};V)$.
We set $\tilde{u}_{a_n}:=u_a(\cdot +a-a_n)\in K_{a_n}$.
Then, $F(a_n;V)\le I(\tilde{u}_{a_n};V)$ holds.
Applying Lemma~\ref{Lebsgue},
we have
\begin{align}\label{Fah}
\limsup_{n\to \infty}F(a_n;V)\le \lim_{n\to \infty}I(\tilde{u}_{a_n};V)=I(u_a;V)=F(a;V).
\end{align}
In particular, we obtain
\[
\sup_{n\in\N} I(u_{a_n};V)=\sup_{n\in\N} F(a_n;V)<\infty.
\]
Then, from Lemma~\ref{newlem}, it follows that
$\{u_{a_n}\}_{n\in\N}$ is bounded in $H^1(\R)$.
Therefore, 
replacing $\{a_n\}_{n\in\N}$ by a subsequence if necessary,
there exists $v_a\in H^1(\R)$ such that
$u_{a_n}$ weakly converges to $v_a$ in $H^1(\R)$ as $n\to \infty$.
It implies $v_a\in K_a$ and
\begin{align*}
I(v_a;V)\le \liminf_{n\to\infty}I(u_{a_n};V).\label{liminfF}
\end{align*}
Thus, it holds that
\[
F(a;V)\le I(v_a;V)
\le \liminf_{n\to\infty}I(u_{a_n};V)=\liminf_{n\to\infty}F(a_n;V)
\le \limsup_{n\to\infty}F(a_n;V)\le F(a;V).
\]
Hence, we obtain $\lim_{n\to\infty}F(a_n;V)=F(a;V)$
and conclude that $F\in C^0(\R)$.
\qed

The following theorem gives a sufficient condition for the trapped mode
by a potential well in terms of the width and depth of the potential well.
\begin{Th}\label{pw}
Let $\alpha,\beta >0$ and $b<c$.
We suppose that $V\in L^\infty(\R)\subset X$
satisfies $V(x)\ge \alpha$ for a.e. $x\in (-\infty,b)\cup (c,\infty)$
and $V(x)=-\beta $ for $x\in (b,c)$.
If $\sqrt{\beta}(c-b)\ge \pi$, then
there exists $a\in [b,c]$ such that $F(a;V)=m(V)$.
\end{Th}
\begin{Rem}
{\rm 
In Theorem~\ref{pw}, $c-b$ represents the width of the potential well,
and $\beta$ is the depth of the potential well.
The condition $\sqrt{\beta}(c-b)\ge \pi$ gives a sufficient condition
for the trapped mode in terms of the width and depth of the potential well.}
\end{Rem}
\noindent
{\it Proof of Theorem~\ref{pw}.}
For $a\in\R$, from Theorems~\ref{fms1} and \ref{fms2}, 
there exists $u_a\in K_a$ such that 
$I(u_a;V)=F(a;V)$ and $u_a(x)>0$ for $x\in \R$.
Let $a\in \R\setminus [b,c]$. 
If $u_a(x)<1$ for $x\in [b,c]$, then,
from Theorem~\ref{fms2},
it has to satisfy $u_a''(x)+\beta u_a(x)=0$
and
\begin{align}\label{sin}
  u_a(x)=C \sin (\sqrt{\beta}x+\theta)\quad (x\in (b,c)),
\end{align}
where
$C\in\R$ and $\theta\in \R$ are some constants. But it is impossible for a function of the
form \eqref{sin} to
satisfy the condition $0<u_a(x)<1$ for $x\in [b,c]$
if $c-b\ge \pi/\sqrt{\beta}$.
Hence, there exists $\tilde{a}\in [b,c]$ such that
$u_a\in K_{\tilde{a}}$ holds.
Since $F(\tilde{a};V)\le I(u_a;V)=F(a;V)$, we conclude that
\begin{align}\label{claim1}
  ~^\forall a\in \R\setminus [b,c],~
  ^\exists \tilde{a}\in [b,c]~~\text{s.t.}~~
  F(\tilde{a};V)\le F(a;V).
\end{align}
Since $F(\cdot;V)\in C^0(\R)$ holds from Theorem~\ref{FC0},
we obtain
\[
\inf_{a\in\R}F(a;V)=\inf_{a\in [b,c]}F(a;V)=\min_{a\in [b,c]}F(a;V).
\]
Therefore, from Theorem~\ref{decomposition}, there exists $a\in [b,c]$ such that $F(a;V)=m(V)$.
\qed

\section*{Acknowledgements}
This work was partially supported by JSPS KAKENHI Grant Nos. 20KK0058,
20H01812, and 20K03675.


\bigskip

\begin{thebibliography}{99}
\bibitem{AKOp1}
V. Apriliani, M. Kimura, and H. Ohtsuka:
Two-step minimization approach to Sobolev-type inequality
with bounded potential in 1D. (arXiv:2311.00708v2)

\bibitem{Bre11} H. Brezis:
Functional Analysis, Sobolev Spaces and Partial Differential Equations.
Springer-Verlag, New York, (2011).

\bibitem{CTZ19}
D. Cassani, C. Tarsi, and J. Zhang:
Bounds for best constants in subcritical Sobolev embeddings.
Nonlinear Analysis, Vol.187 (2019), 438-449.

\bibitem{CM16}
A. Cianchi and V. Maz'ya:
Sobolev inequalities in arbitrary domains.
Advances in Mathematics, Vol.293 (2016), 644-696.

\bibitem{Kam08}
Y. Kametaka, H. Yamagishi, K. Watanabe, A. Nagai, K. Takemura, and M. Arai:
The best constant of Sobolev inequality which corresponds
to Schr\"odinger operator with Dirac delta potential.
Scientiae Mathematicae Japonicae Online, Vol.e-2008 (2008), 541-555.


\bibitem{WKNTY08}
K. Watanabe, Y. Kametaka, A. Nagai, K. Takemura, and H. Yamagishi:
The best constant of Sobolev inequality on a bounded interval.
J. Math. Anal. Appl., Vol.340 (2008), 699-706.

\bibitem{Will07}
M. Willem: Functional Analysis
Fundamentals and Applications, Second Edition.
Birkh\"{a}user New York, NY (2007).

\end{thebibliography}
\end{document}